\documentclass{article}
\usepackage{amssymb,amsmath,theorem,euscript}

\input{epsf}

\newcounter{sec}

\def\sm{\smallskip}

\newcounter{punct}[sec]

\def\punct{\refstepcounter{punct}{\arabic{sec}.\arabic{punct}.  }} 

\def\COUNTERS{\addtocounter{sec}{1}
              \setcounter{punct}{0}
          \setcounter{equation}{0}
          \setcounter{theorem}{0}
                  }

\newtheorem{theorem}{Theorem}[sec]

\newtheorem{lemma}[theorem]{Lemma}

\begin{document}

 \def\ov{\overline}
\def\wt{\widetilde}
 \newcommand{\rk}{\mathop {\mathrm {rk}}\nolimits}
\newcommand{\Aut}{\mathop {\mathrm {Aut}}\nolimits}
\newcommand{\Out}{\mathop {\mathrm {Out}}\nolimits}
 \newcommand{\tr}{\mathop {\mathrm {tr}}\nolimits}
  \newcommand{\diag}{\mathop {\mathrm {diag}}\nolimits}
  \newcommand{\supp}{\mathop {\mathrm {supp}}\nolimits}
  \newcommand{\indef}{\mathop {\mathrm {indef}}\nolimits}
  \newcommand{\dom}{\mathop {\mathrm {dom}}\nolimits}
  \newcommand{\im}{\mathop {\mathrm {im}}\nolimits}
 
\renewcommand{\Re}{\mathop {\mathrm {Re}}\nolimits}
\renewcommand{\Im}{\mathop {\mathrm {Im}}\nolimits}

\def\Br{\mathrm {Br}}

\def\SL{\mathrm {SL}}
\def\Diag{\mathrm {Diag}}
\def\SU{\mathrm {SU}}
\def\GL{\mathrm {GL}}
\def\U{\mathrm U}
\def\OO{\mathrm O}
 \def\Sp{\mathrm {Sp}}
 \def\SO{\mathrm {SO}}
\def\SOS{\mathrm {SO}^*}
 \def\Diff{\mathrm{Diff}}
 \def\Vect{\mathfrak{Vect}}
\def\PGL{\mathrm {PGL}}
\def\PU{\mathrm {PU}}
\def\PSL{\mathrm {PSL}}
\def\Symp{\mathrm{Symp}}
\def\End{\mathrm{End}}
\def\Mor{\mathrm{Mor}}
\def\Aut{\mathrm{Aut}}
 \def\PB{\mathrm{PB}}
 \def\cA{\mathcal A}
\def\cB{\mathcal B}
\def\cC{\mathcal C}
\def\cD{\mathcal D}
\def\cE{\mathcal E}
\def\cF{\mathcal F}
\def\cG{\mathcal G}
\def\cH{\mathcal H}
\def\cJ{\mathcal J}
\def\cI{\mathcal I}
\def\cK{\mathcal K}
 \def\cL{\mathcal L}
\def\cM{\mathcal M}
\def\cN{\mathcal N}
 \def\cO{\mathcal O}
\def\cP{\mathcal P}
\def\cQ{\mathcal Q}
\def\cR{\mathcal R}
\def\cS{\mathcal S}
\def\cT{\mathcal T}
\def\cU{\mathcal U}
\def\cV{\mathcal V}
 \def\cW{\mathcal W}
\def\cX{\mathcal X}
 \def\cY{\mathcal Y}
 \def\cZ{\mathcal Z}
\def\0{{\ov 0}}
 \def\1{{\ov 1}}
 \def\frA{\mathfrak A}
 \def\frB{\mathfrak B}
\def\frC{\mathfrak C}
\def\frD{\mathfrak D}
\def\frE{\mathfrak E}
\def\frF{\mathfrak F}
\def\frG{\mathfrak G}
\def\frH{\mathfrak H}
\def\frI{\mathfrak I}
 \def\frJ{\mathfrak J}
 \def\frK{\mathfrak K}
 \def\frL{\mathfrak L}
\def\frM{\mathfrak M}
 \def\frN{\mathfrak N} \def\frO{\mathfrak O} \def\frP{\mathfrak P} \def\frQ{\mathfrak Q} \def\frR{\mathfrak R}
 \def\frS{\mathfrak S} \def\frT{\mathfrak T} \def\frU{\mathfrak U} \def\frV{\mathfrak V} \def\frW{\mathfrak W}
 \def\frX{\mathfrak X} \def\frY{\mathfrak Y} \def\frZ{\mathfrak Z} \def\fra{\mathfrak a} \def\frb{\mathfrak b}
 \def\frc{\mathfrak c} \def\frd{\mathfrak d} \def\fre{\mathfrak e} \def\frf{\mathfrak f} \def\frg{\mathfrak g}
 \def\frh{\mathfrak h} \def\fri{\mathfrak i} \def\frj{\mathfrak j} \def\frk{\mathfrak k} \def\frl{\mathfrak l}
 \def\frm{\mathfrak m} \def\frn{\mathfrak n} \def\fro{\mathfrak o} \def\frp{\mathfrak p} \def\frq{\mathfrak q}
 \def\frr{\mathfrak r} \def\frs{\mathfrak s} \def\frt{\mathfrak t} \def\fru{\mathfrak u} \def\frv{\mathfrak v}
 \def\frw{\mathfrak w} \def\frx{\mathfrak x} \def\fry{\mathfrak y} \def\frz{\mathfrak z} \def\frsp{\mathfrak{sp}}
 \def\bfa{\mathbf a} \def\bfb{\mathbf b} \def\bfc{\mathbf c} \def\bfd{\mathbf d} \def\bfe{\mathbf e} \def\bff{\mathbf f}
 \def\bfg{\mathbf g} \def\bfh{\mathbf h} \def\bfi{\mathbf i} \def\bfj{\mathbf j} \def\bfk{\mathbf k} \def\bfl{\mathbf l}
 \def\bfm{\mathbf m} \def\bfn{\mathbf n} \def\bfo{\mathbf o} \def\bfp{\mathbf p} \def\bfq{\mathbf q} \def\bfr{\mathbf r}
 \def\bfs{\mathbf s} \def\bft{\mathbf t} \def\bfu{\mathbf u} \def\bfv{\mathbf v} \def\bfw{\mathbf w} \def\bfx{\mathbf x}
 \def\bfy{\mathbf y} \def\bfz{\mathbf z} \def\bfA{\mathbf A} \def\bfB{\mathbf B} \def\bfC{\mathbf C} \def\bfD{\mathbf D}
 \def\bfE{\mathbf E} \def\bfF{\mathbf F} \def\bfG{\mathbf G} \def\bfH{\mathbf H} \def\bfI{\mathbf I} \def\bfJ{\mathbf J}
 \def\bfK{\mathbf K} \def\bfL{\mathbf L} \def\bfM{\mathbf M} \def\bfN{\mathbf N} \def\bfO{\mathbf O} \def\bfP{\mathbf P}
 \def\bfQ{\mathbf Q} \def\bfR{\mathbf R} \def\bfS{\mathbf S} \def\bfT{\mathbf T} \def\bfU{\mathbf U} \def\bfV{\mathbf V}
 \def\bfW{\mathbf W} \def\bfX{\mathbf X} \def\bfY{\mathbf Y} \def\bfZ{\mathbf Z} \def\bfw{\mathbf w}
 \def\R {{\mathbb R }} \def\C {{\mathbb C }} \def\Z{{\mathbb Z}} \def\H{{\mathbb H}} \def\K{{\mathbb K}}
 \def\N{{\mathbb N}} \def\Q{{\mathbb Q}} \def\A{{\mathbb A}} \def\T{\mathbb T} \def\P{\mathbb P} \def\G{\mathbb G}
 \def\bbA{\mathbb A} \def\bbB{\mathbb B} \def\bbD{\mathbb D} \def\bbE{\mathbb E} \def\bbF{\mathbb F} \def\bbG{\mathbb G}
 \def\bbI{\mathbb I} \def\bbJ{\mathbb J} \def\bbK{\mathbb K} \def\bbL{\mathbb L} \def\bbM{\mathbb M} \def\bbN{\mathbb N} \def\bbO{\mathbb O}
 \def\bbP{\mathbb P} \def\bbQ{\mathbb Q} \def\bbS{\mathbb S} \def\bbT{\mathbb T} \def\bbU{\mathbb U} \def\bbV{\mathbb V}
 \def\bbW{\mathbb W} \def\bbX{\mathbb X} \def\bbY{\mathbb Y} \def\kappa{\varkappa} \def\epsilon{\varepsilon}
 \def\phi{\varphi} \def\le{\leqslant} \def\ge{\geqslant}

\def\UU{\bbU}
\def\Mat{\mathrm{Mat}}
\def\tto{\rightrightarrows}

\def\Gr{\mathrm{Gr}}

\def\graph{\mathrm{graph}}

\def\O{\mathrm{O}}

\def\la{\langle}
\def\ra{\rangle}

\def\B{\mathrm B}
\def\Int{\mathrm{Int}}
\def\LGr{\mathrm{LGr}}


\def\I{\mathbb I}
\def\M{\mathbb M}
\def\T{\mathbb T}
\def\S{\mathrm S}

\def\Lat{\mathrm{Lat}}
\def\LLat{\mathrm{LLat}} 

\def\MAT{\mathrm{MAT}}
\def\Mar{\mathrm{Mar}}

\def\Bd{\mathrm{Bd}}
\def\We{\mathrm{We}}
\def\Heis{\mathrm{Heis}}
\def\Pol{\mathrm{Pol}}
\def\Ams{\mathrm{Ams}}
\def\Herm{\mathrm{Herm}}

\def\bbot{{\bot\!\!\!\bot}}

\begin{center}
 \Large \bf
 
 Whishart--Pickrell distributions and closures of group actions 
 
 \large \sc
 
 \bigskip
 
 Yu.A.Neretin%
 \footnote{Supported by grant FWF, P28421, P25142}
 
\end{center}

{\small Consider probabilistic distributions on the space 
of infinite Hermitian matrices $\Herm(\infty)$
invariant with respect
to the unitary group $\U(\infty)$. We describe the closure
of $\U(\infty)$ in the space of spreading maps (polymorphisms) of $\Herm(\infty)$,
this closure is a semigroup isomorphic to the semigroup of all contractive operators.}

\section{The statement}

\COUNTERS

{\bf\punct Notation.}
Denote by $\Herm_\infty$ the space of all infinite Hermitian matrices.
By $\Herm_\infty^0$ we denote the space of all infinite Hermitian matrices having
a finite number of non-zero matrix elements.

By $\U(\infty)$ we denote the group  of  infinite unitary matrices
$g$ such that 
 $g-1$ has only finite number of nonzero matrix elements.
 By $\ov\U(\infty)$ we denote the complete  unitary group in $\ell^2$,
 we equip it with the weak operator topology. 
 
 By $\cB(\infty)$ we denote the semigroup of all linear operators
 in $\ell_2$ with norm $\le 1$.
 
 \sm

{\bf\punct  Whishart--Pickrell distributions.}
For any probabilistic measure $\mu$ on $\Herm_\infty$ we assign the
characteristic function on $\Herm_\infty^0$ by
$$
\chi(\mu|A)=\int_{\Herm_\infty} e^{i\tr AX} \,d\mu(X).
$$
Obviously, such a function determines a measure in a unique way.

Consider the group $\U(\infty)$ of infinite unitary matrices
$g$ such that 
 $g-1$ has only finite number of nonzero matrix elements.
This group acts on $\Herm_\infty$ by conjugations,
$$
U:\, X\mapsto U^{-1} X U.
$$

There is the following theorem of Pickrell \cite{Pick} (see also
another proof with additional details in \cite{OV}) 
in the spirit of the de Finetti theorem:

\sm

{\it Any $\U(\infty)$-invariant measure on $\Herm_\infty$ can be decomposed into ergodic measures.
An ergodic $\U(\infty)$-invariant measure has a characteristic function of the form 
\begin{equation}
 \chi_{\gamma_1,\gamma_2,\lambda}(A)
 =
 e^{-\frac{\gamma_1} 2\tr A^2+i\gamma_2 \tr A}
\prod_{k=1}^\infty \Bigl( \det {\frac{e^{-i \lambda_k A}}{1-i \lambda_k A}} \Bigr),
\label{eq:characteristic}
\end{equation}
where $\gamma_1\ge 0$, $\gamma_2$, and $\lambda_1$, $\lambda_2$, \dots
are real numbers, and $\sum\lambda_k^2<\infty$.
}

\sm

Denote this measure by $\mu_{\gamma_1,\gamma_2,\lambda}$.

\sm

A characteristic function is a product and the corresponding measure can be decomposed as an (infinite)
convolution. The factor $ e^{-\frac{\gamma_1} 2\tr A^2+i\gamma_2 \tr A}$ corresponds to a Gaussian measure
on $\Herm_\infty$. Let us explain a meaning of factors $\det(1-i\lambda_k A)^{-1}$. Consider the space
$\C$ equipped with the Gaussian measure $\pi^{-1} e^{-|u|^2}\, d\Re u\,d\Im u$. Concider the space $\C^\infty$
equipped with a product-measure $\nu$.
Concider the map $\C^\infty\to \Herm_\infty$
given by
$$
u\mapsto \lambda_k u^*u.
$$
Concider the image $\mu$ of $\nu$ under this map. The characteristic function of $\mu$ is
\begin{multline*}
\int_{\Herm_\infty} e^{i\tr AX}d\mu(X)=\int_{\C^\infty} e^{i\tr \lambda_k Au^*u}\,d\nu(u)=
\int_{\C^\infty} e^{i \lambda_k uAu^*}\,d\nu(u)=\\= \det(1-i\lambda_k A)^{-1}.
\end{multline*}

If $\sum |\lambda_k|<\infty$, then we can transform the expression (\ref{eq:characteristic})
to
\begin{equation*}
 \chi_{\gamma_1,\gamma_2,\lambda}(A)
 =
 e^{-\frac{\gamma_1} 2\tr A^2+i(\gamma_2-\sum \lambda_k) \tr A}
\prod_{k=1}^\infty  \det(1-i \lambda_k A)^{-1}.
\end{equation*}
If the series $\sum \lambda_k$ diverges, we get a divergent series in the exponent and a divergent product.

\sm

{\bf \punct Polymorphisms.}  See \cite{Kre}, \cite{Ver}, \cite{Ner-book}, Section VIII.4.
Consider a Lebesgue measure space $M$ with a non-atomic probabilistic measure $\mu$.
Denote by $\Ams(M)$ the group of measure preserving bijective a.s.
transformations of $M$.

A {\it polymorphism} of $M$ is a measure $\kappa$ on $M\times M$, whose pushforwards
to $M$ with respect to the both projections $M\times M\to M$ coincide with $\mu$. 
Denote by $\Pol(M)$ the set of all polymorphisms of $M$. 
We say that a  sequence $\pi_j\in\Pol(M)$ converges to $\pi$ if for any measurable sets
$A$, $B\subset M$ we have convergence $\pi_j(A\times B)\to\pi(A\times B)$).
The space $\Pol(M)$ is compact and the group $\Ams(M)$ is dense
in $\Pol(M)$.

Polymorphisms can be regarded as maps, spreading points of $M$ to probabilistic
measures on $M$. Namely, for $\pi\in\Pol(M)$ consider a system of  conditional measures
$\pi_m$ on sets $m\times M\subset M\times M$, where $m$ ranges in $M$.
We declare that the 'map' $\pi$ send each point $m$ to the measure $\pi_m$.
If $\pi$, $\kappa\in \Pol(M)$, then the product $\rho=\kappa\circ \pi$
is defined from the condition
$$
\rho_m=\int_M \kappa_n \,d\pi_m(n).
$$
We get a semigroup with separately continuous product.

For any $g\in\Ams(M)$ consider the map from $M$ to $M\times M$ defined by $m\mapsto (m,g(m))$
and take the pushforward of the measure $\mu$.
In this way, we get  a polymorphism supported by the graph of $g$.
The group $\Ams(M)$ is dense
in $\Pol(M)$.

A  {\it Markov operator} $R$ in $L^2(M)$ is a bounded operator satisfying 
3 properties: 

\sm

-- for any function $f\ge 0$ we have $Rf\ge 0$;

\sm

-- $R\cdot 1=1$, $R^*\cdot 1=1$.

\sm

Recall that automatically $\|R\|= 1$.
There is a one-to-one correspondence between the set of Markov operators $\Mar(M)$
and $\Pol(M)$. Namely, let $R$ be a Markov operator, then we define a polymorphism
$\pi$ by
$$
\pi(A\times B)=\la R I_A,I_B\ra_{L^2(M)},
$$
where $A$, $B\subset M$ are measurable sets, and $I_A$, $I_B$ are their
indicator functions.
The weak convergence in $\Mar(M)$
corresponds to the  convergence in $\Pol(M)$,
the product of Markov operators corresponds
to the product of polymorphisms.

\sm

{\bf\punct Closures of actions.} Let a group $G$ act on
$M$ by measure-preserving transformations, i.e., we have 
a homomorphism $G\to\Ams(M)$, to be definite assume that 
this is an embedding. Then the closure of $G$ in $\Pol(M)$ is a compact semigroup
$\Delta\supset G$. A description of the closure is not too interesting for connected
Lie groups (for instance, for semisimple linear Lie groups we get a one-point compactification,
this follows from  \cite{HM}, Theory 5.3). However, a description
of such closure is interesting  for infinite-dimensional groups.

The first result of this type was obtained by Nelson
\cite{Nel}, 1973. He showed that the standard action of the infinite-dimensional
orthogonal group on the space with Gaussian measure admits an extension to an action of the semigroup
of all contractive operators by polymorphisms and obtained formulas for such measures.
Now, it is known a big collection of actions infinite-dimensional groups 
on measure spaces, however the closures are evaluated in few cases \cite{Ner-gauss}, \cite{Ner-poisson}.
Here we describe a new (relatively simple) example.

\sm

{\bf\punct The statement.}
For any polymorphism
$\pi$ on $\Herm_\infty$ we write a characteristic function
on $\Herm_\infty^0\times\Herm_\infty^0$ by
$$
F(\pi|A,B):=\int_{\Herm_\infty\times\Herm_\infty}
e^{i\tr AX+BY}\,d\pi(X,Y).
$$

Denote by $\cB(\infty)$ the semigroup of all operators in $\ell_2$ with norm $\le 1$.
Our purpose  is the following statement

\sm

\begin{theorem}
For any ergodic measure on $\Herm_\infty$, the closure of $\U(\infty)$ in $\Pol(\Herm_\infty)$
is isomorphic to $\cB(\infty)$.  If $S\in \cB(\infty)$,
then the characteristic function of the corresponding polymorphism $\pi_S$
is given by
\begin{multline}
 F(\pi_S|A,B)= \exp\Bigl\{-\frac {\gamma_1}2 \tr
 \Bigl[\begin{pmatrix}A&0\\0&B \end{pmatrix} 
 \begin{pmatrix} 1&S\\ S^*&1  \end{pmatrix}\Bigr]^2+i\gamma_2(\tr A+\tr B)
   \Bigr\}
   \times\\\times
   \prod_k \frac{e^{-i\lambda_k(\tr A+\tr B)}}{\det\Bigl[1-i\lambda_k 
 \begin{pmatrix}A&0\\0&B \end{pmatrix} 
 \begin{pmatrix} 1&S\\ S^*&1  \end{pmatrix}
   \Bigr]}
   \,.
   \label{eq:answer}
\end{multline}
\end{theorem}

 Since the characteristic function is a product,
the measure $\pi_S$ can be decomposed as a convolution 
of measures. The exponential factor corresponds to a Gaussian measure,
let us explain meaning of factors in the product.
Concider a measure $\nu_S$ on $\C^\infty\times \C^\infty$ defined in the following way in terms of its characteristix function
$$
\int_{\C^\infty\times \C^\infty} e^{i\Re u \ov z+i\Re v\ov z} d\nu_S:=
\exp\Bigl\{-\frac 12 \begin{pmatrix} u&v\end{pmatrix}  \begin{pmatrix} 1&S\\ S^*&1  \end{pmatrix}
\begin{pmatrix} u^*\\v^*\end{pmatrix} 
\Bigr\}.
$$
Consider a map $\C^\infty\times \C^\infty\to \Herm_\infty\times \Herm_\infty$
given by 
$$
(u,v)\mapsto (\lambda u^* u,\lambda v^*v).
$$
Then the image of $\nu_S$ under this map is a measure whose characteristic function is
$$
\det\Bigl[1-i\lambda_k 
 \begin{pmatrix}A&0\\0&B \end{pmatrix} 
 \begin{pmatrix} 1&S\\ S^*&1  \end{pmatrix}
   \Bigr]^{-1}.
$$

\section{Proof}

\COUNTERS

{\bf\punct A priori remarks.}

\begin{theorem}
\label{th:2}
 Let the complete unitary group $\ov \U(\infty)$ act by measure preserving transformations%
 \footnote{Such an action can not be point-wise;
 the transformations are defined a.s., and products also are defined a.s.,
 see \cite{Gla}.}
 on a Lebesgue space $M$ with a probability measure. Then the closure of 
 $\ov \U(\infty)$ in $\Pol(M)$ is isomorphic to $\cB(\infty)$.
\end{theorem}

{\sc Proof.} Let $\rho$ be a unitary representation of $\ov\U(\infty)$
in a Hilbert space $H$.
 The closure of the group $\rho(\ov\U(\infty))$
in the group of all bounded operators in $H$ with respect to the weak operator topology
is a semigroup isomorphic to $\cB(\infty)$ (this follows from Kirillov--Olshanski
classification of representations of  $\ov\U(\infty)$, see \cite{Olsh}, Theorem 1.2).

We apply this to the action of  $\ov\U(\infty)$ in $L^2(M)$.
The group $\ov\U(\infty)$ acts by Markov operators, and weak limits of Markov operators
are Markov operators. therefore the semigroup
$\cB(\infty)$ also acts by Markov operators.\hfill $\square$

\sm

\begin{lemma}
 Any ergodic action of $\U(\infty)$ on $\Herm_\infty$ admits a continuous extension
 to an action of $\ov\U(\infty)$.
\end{lemma}

{\sc Proof.} According \cite{OV}, Corollary 2.14, the representation
of $\U(\infty)$ in the space $L^2(\Herm_\infty,\mu)$ has a continuous extension to $\ov\U(\infty)$.
By the  continuity arguments, the group  $\ov\U(\infty)$ acts by Markov operators $\rho(g)$.
We have $\|\rho(g)\|\le 1$, $\|\rho(g)^{-1}\|\le 1$. Therefore $\rho(g)$ is a unitary operator.
Hence it corresponds to a measure preserving transformation.
\hfill $\square$

\sm

{\bf\punct Calculation.} 
Consider a measure $\mu$ with a characteristic function 
(\ref{eq:characteristic}). 
For $g\in\U(\infty)$ consider the corresponding polymorphism $\pi_g$
of $\Herm_\infty$.
The characteristic function of $\pi_g$ is
\begin{multline}
F(\pi_g|A,B)=\exp\Bigl\{-\frac{\gamma_1}2 \tr (A+ UBU^{-1})^2+i\gamma_2(\tr A+\tr U BU^{-1}) \Bigr\}
\times\\ \times
\prod_{k=1}^\infty\frac{e^{-i\lambda_k \tr (A+UBU^{-1})}}
{\det\bigl[1-i\lambda_k (A+UBU^{-1}) \bigr]}.
\label{eq:FAB}
\end{multline}

By Theorem \ref{th:2}, for any $R\in \cB(\infty)$,
we have a polymorphism $\pi_R$ in the space $L^2(\Herm_\infty,\mu)$.
If a sequence $R_j\in\cB(\infty)$ weakly converges to $R$, then we have the weak convergence of
the corresponding polymorphisms, $\pi_{R_j}\to \pi_R$. This is equivalent to a point-wise
convergence of characteristic functions. Now for $R\in \cB(\infty)$ and  a sequence
$g_j\in\U(\infty)$ weakly converging to $R$ we can find 
$F(\pi_R|A,B)$ as a point-wise limit of 
 $F(\pi_{g_j}|A,B)$.

Let $S$ be a finitary operator with norm $\le 1$, let actually
$S$ be $(\alpha+\infty)$-block matrix of the form
$$
 S=
\begin{pmatrix}
u&0\\0& 0
\end{pmatrix}
.
$$
Then we can build $u$ as a block of a unitary $(\alpha+\alpha)$-block matrix 
$\begin{pmatrix}
  u&v\\ w&z                                                        
\end{pmatrix}
$ (see, e.g., \cite{Ner-book}, Theorem VIII.3.2).
Choose $U=U_m$ as a block unitary $\alpha+m+m+\alpha+\infty$-matrix of the form
$$
U_m=\begin{pmatrix}
   u&0&0&v&0\\
   0&0&1&0&0\\
   0&1&0&0&0\\
   w&0&0&z&0\\
   0&0&0&0&1
  \end{pmatrix}.
$$
Clearly, we have a weak convergence
$$
U_m\to S,
$$
 we wish to watch the convergence of characteristic functions
 $F(\pi_{U_m}|A,B)$.
In fact, we will show that this sequence is eventually constant
for any fixed $A$, $B$.

Fix $A$, $B\in\Herm_\infty^0$. Let actually $A$, $B\in\Herm_{\alpha+\beta}$.
Let $m$ be sufficiently large (in fact, $m\ge\beta$). We represent the matrices $A$, $B$ as block matrices
of size $\alpha+\beta+(m-\beta)+\beta+(m-\beta)+\infty$:
$$
A=\begin{pmatrix}
   a_{11}&a_{12}&0&\dots&0\\
    a_{21}&a_{22}&0&\dots&0\\
    0&0&0&\dots&0\\
    \vdots&\vdots&\vdots&\ddots&\vdots\\
    0&0&0&\dots&0\\
  \end{pmatrix}
,\qquad 
B=\begin{pmatrix}
   b_{11}&b_{12}&0&\dots&0\\
    b_{21}&b_{22}&0&\dots&0\\
    0&0&0&\dots&0\\
    \vdots&\vdots&\vdots&\ddots&\cdots\\
    0&0&0&\dots&0\\
  \end{pmatrix}.
$$
We wish to evaluate 
$$
\det\bigl[1-i\lambda_k (A+U_mBU_m^{-1}) \bigr] \qquad\text{and}\qquad
\tr (A+U_mBU_m^{-1})^2.
$$
A straightforward calculation gives
$$
A+U_mBU_m^{-1}=\begin{pmatrix}
               a_{11}+ ub_{11} u^*& a_{12}&0& ub_{12}& ub_{11}w^*&0&0\\            
               a_{21}&a_{22}&0&0& 0&0&0\\
                 0&0&0&0&0&0&0\\
               b_{21}u^*&0&0&b_{22}&b_{21}w^*&0&0\\
               wb_{11} u^*&0&0&wb_{12}&wb_{11}w^*&0&0\\
               0&0&0&0&0&0&0\\
               0&0&0&0&0&0&0
               \end{pmatrix}.
$$
Clearly, we can remove zero columns and zero rows from this matrix. Formally,
$$
\det\bigl(1-i\lambda_k (A+U_mBU_m^{-1})\bigr)=\det(1-i\lambda_k H), 
$$
where
$$
H=
\begin{pmatrix}
               a_{11}+ ub_{11} u^*& a_{12}& ub_{12}& ub_{11}w^*\\               
               a_{21}&a_{22}&0& 0\\
               b_{21}u^*&0&b_{22}&b_{21}w^*\\
               wb_{11} u^*&0&wb_{12}&wb_{11}w^*\\             
               \end{pmatrix}.
$$
Denote by $\Delta$ the block diagonal matrix with blocks $1$, $1$, $1$, $w$.
Represent $H$ as $H=\Delta Z$ (where the expression for $Z$ is clear).
We apply the formula
$$
\det(1-i\lambda_k \Delta Z)=\det(1-i\lambda_k Z\Delta)
.
$$
Denote $H':=Z\Delta $,
\begin{equation}
H'=
\begin{pmatrix}
               a_{11}+ ub_{11} u^*& a_{12}& ub_{12}& ub_{11}w^*w\\
               a_{21}&a_{22}& 0&0\\
               b_{21}u^*&0&b_{22}&b_{21}w^*w\\
               b_{11} u^*&0&b_{12}&b_{11}w^*w\\
               \end{pmatrix}
               .
               \label{eq:step1}
\end{equation}
Since the matrix $\begin{pmatrix}
  u&v\\ w&z                                                        
\end{pmatrix}$ is unitary, we have
 $w^*w=1-uu^*$. We substitute this to the 4-th column of (\ref{eq:step1}),
 keeping a result  in mind we continue our calculations.
 Denote
 $$
 T=\begin{pmatrix}
    1&0&0&u\\
    0&1&0&0\\
    0&0&1&0\\
    0&0&0&1\\
   \end{pmatrix}
 .$$
 Then
 \begin{multline*}
  H'=
  \begin{pmatrix}
               a_{11}+ ub_{11} u^*& a_{12}& a_{11}u+ ub_{12}& ub_{11}\\
               a_{21}&a_{22}& 0&a_{21}u\\
               b_{21}u^*&0&b_{22}&b_{21}\\
               b_{11} u^*&0&b_{12}&b_{11}\\
               \end{pmatrix} T^{-1}=\\=
              T \begin{pmatrix}
 a_{11}& a_{12}&0& a_{11} u\\
 a_{21} &a_{22}&0& a_{21}u\\
 b_{21} u^*&0&b_{22}&b_{21}\\
  b_{11} u^*&0&b_{12}&b_{11}
\end{pmatrix} T^{-1}.
 \end{multline*}

 Hence 
 \begin{multline}
 \det\left(1-i\lambda_k H'\right)=
 \det\left[1-i\lambda_k
\begin{pmatrix}
 a_{11}& a_{12}&0& a_{11} u\\
 a_{21} &a_{22}&0& a_{21}u\\
 b_{21} u^*&0&b_{22}&b_{21}\\
  b_{11} u^*&0&b_{12}&b_{11}
\end{pmatrix}\right]
=\\=
\det\left[1-i\lambda 
\begin{pmatrix}
 a_{11}& a_{12}& a_{11} u&0\\
 a_{21} &a_{22}&a_{21}u&0\\
 b_{11} u^*&0&b_{11}&b_{12}\\
 b_{21} u^*&0&b_{21}&b_{22}\end{pmatrix}\right]
=\\=
\det\left[1-i\lambda \begin{pmatrix}A&0\\ 0&B \end{pmatrix}
\begin{pmatrix}
 1&S\\ S^*&1
\end{pmatrix}
\right]
\end{multline}
and we get the desired expression.

Next,
$$
\tr (A+U_mBU_m^{-1})^2=\tr A^2+\tr B^2+2\tr A U_mBU_m^{-1}.
$$
Multiplying  matrices we observe that 
$A U_mBU_m^{-1}$ has a unique nonzero diagonal block, $a_{11}u b_{11} u^*$.
Thus, 
$$
\tr A U_mBU_m^{-1}=\tr a_{11}u\, b_{11} u^*= \tr A S B S^*,
$$
and this implies 
$$
\tr (A+U_mBU_m^{-1})^2=
\tr
\left[
\begin{pmatrix}A&0\\ 0&B \end{pmatrix}
\begin{pmatrix}
 1&S\\ S^*&1
\end{pmatrix}\right]^2.
$$

Thus, for $m\ge \beta$ the value of
$F(\pi_{U_m}|A,B)$ is given by formula (\ref{eq:answer}),

So, the theorem holds for finitary matrices $S$.
Consider an arbitrary  $S\in \cB(\infty)$.
Denote by $S[m]$ the left upper corner of $S$.
Denote
$$
S_m:=\begin{pmatrix}
      S[m]&0\\0&0
     \end{pmatrix}.
$$
We have a weak convergence $S_m\to S$. On the other hand,
we have the pointwise convergence
$$
F(\pi_{S_m}|A,B)\to F(\pi_{S}|A,B)
$$
(in fact, this sequence is eventually constant for
any fixed $A$, $B$).

This completes the proof of the theorem.

\noindent
\tt Math.Dept., University of Vienna,
 \\
 Oskar-Morgenstern-Platz 1, 1090 Wien;
 \\
\& Institute for Theoretical and Experimental Physics (Moscow);
\\
\& Mech.Math.Dept., Moscow State University;
\\
Institute for information transmission problems (Moscow);
\\
e-mail: neretin(at) mccme.ru
\\
URL:www.mat.univie.ac.at/$\sim$neretin

\end{document}